\documentclass[10pt, twoside]{article}
\usepackage{amssymb,amsmath}
\usepackage[amsmath,thmmarks,thref]{ntheorem}
\usepackage{graphics}
\usepackage{enumerate}
\usepackage{mathrsfs}
\usepackage{dsfont}
\title{Filtered Colimit Preserving Functors on Models of a Regular Theory}
\author{Henrik Forssell\thanks{Dept.\ of Math., Stockholm University}
}

\input{diagxy}
\xyoption{curve}
\newbox\anglebox 
\setbox\anglebox=\hbox{\xy \POS(75,0)\ar@{-} (0,0) \ar@{-} (75,75)\endxy}
 
\newbox\angleboxr 
\setbox\angleboxr=\hbox{\xy \POS(0,0)\ar@{-} (0,75) \ar@{-} (75,0)\endxy}
 
\newbox\sanglebox 
\setbox\sanglebox=\hbox{\xy \POS(50,0)\ar@{-} (0,0) \ar@{-} (50,50)\endxy}
 \def\spbangle{\copy\sanglebox}
\newbox\sangleboxr 
\setbox\sangleboxr=\hbox{\xy \POS(0,0)\ar@{-} (0,50) \ar@{-} (50,0)\endxy}
 
\newbox\sangleboxf 
\setbox\sangleboxf=\hbox{\xy \POS(50,50)\ar@{-} (50,0) \ar@{-} (0,50)\endxy}
 
\newbox\angleboxf 
\setbox\angleboxf=\hbox{\xy \POS(75,75)\ar@{-} (75,0) \ar@{-} (0,75)\endxy}
 
\newbox\sangleboxfr 
\setbox\sangleboxfr=\hbox{\xy \POS(0,50)\ar@{-} (50,50) \ar@{-} (0,0)\endxy}
 
\newbox\angleboxfr 
\setbox\angleboxfr=\hbox{\xy \POS(0,75)\ar@{-} (75,75) \ar@{-} (0,0)\endxy}

\newdir{|>}{!/4.7pt/\dir{|}
        *:(1,-.2)\dir^{>}
        *:(1,+.2)\dir_{>}}

\def\embedd{\to/^{ (}->/}
\newcommand{\colimit}{\underrightarrow{\lim}}

\newcommand{\pair}[1]{\ensuremath{\langle {#1} \rangle}}
\newcommand{\homset}[3]{\ensuremath{\operatorname{Hom}_{#1}\!\left({#2},{#3}\right)}}

\newcommand{\mb}{\ensuremath{\mathbin}}

\newcommand{\cat}[1]{\ensuremath{\mathcal{#1}}}

\newcommand{\synt}[2]{\ensuremath{\mathcal{#1}_{\mathbb{#2}}}}
\newcommand{\alg}[1]{\ensuremath{\mathbf{#1}}}

\newcommand{\modcat}[1]{\ensuremath{\mathrm{Mod}_{#1}}}

\newcommand{\modin}[2]{\ensuremath{\mathrm{Mod}_{#1}(#2)}}
\newcommand{\topo}[1]{\ensuremath{\mathscr{#1}}}
\newcommand{\Sets}{\ensuremath{\mathbf{Set}}}

\newcommand{\Sh}[1]{\protect\ensuremath{\operatorname{Sh}\!\left(#1\right)}}

\newcommand{\sem}[1]{\ensuremath{[\![{#1}]\!]}}
\newcommand{\csem}[2]{\ensuremath{[\![{#1}\mb|{#2}]\!]}}

\newcommand{\fins}[1]{\exists {#1}\mathpunct .}

\newcommand{\theory}{\ensuremath{\mathbb{T}}}

\newcommand{\cterm}[2]{\ensuremath{\left \{ {#1}\ \; \vrule \; \ {#2}\right \}}}
\newcommand{\syntob}[2]{\ensuremath{[{#1}\;|\;{#2}]}}
\newcommand{\classtop}{\ensuremath{\mathbf{Set}[\theory]}}

\newcommand{\Eqsheav}[2]{\protect\ensuremath{\operatorname{Sh}_{#1}(#2)}}
\newcommand{\sox}[2]{\ensuremath{\csem{#1}{#2}_{M_{\theory}}}}
\newcommand{\bopen}[1]{\ensuremath{\langle\!\! \langle {#1} \rangle\!\! \rangle}}
\makeindex 
\begin{document}
%
%
\maketitle
\begin{abstract}
This note recalls the representation of  regular theories \theory\ in terms of set-valued functors on models given by Makkai in \cite{makkai:90}, and explicitly states the representation theorem for the classifying topos \classtop\ in terms of filtered colimit preserving functors which can be extrapolated from the results of that paper. That representation of \classtop\ is then compared with topological representations in the style of Butz and Moerdijk \cite{butz:98b} by showing that for a certain natural topology on the space of models, preserving filtered colimits is the same thing as being `continuous' in the sense of being an equivariant sheaf. By using a slight variation of the topology originally presented in \emph{op.\ cit.}, we obtain from this comparison a representation of \classtop\ in terms of a topological category of models and homomorphisms, where the restricted topological groupoid of models and isomorphisms classifies a different (non-regular) theory.

\end{abstract}
%
%
\section{Introduction}
\label{introduction}

By \emph{regular logic} is understood the fragment of first-order logic formed by sequents $\phi\vdash_{\alg{x}}\psi$ of formulas $\phi$, $\psi$ formed with the connectives $\top$, $\wedge$, and $\exists$ (see e.g.\ \cite{elephant1}). An example of a regular theory---that is, one axiomatized by or consisting of regular sequents---is the theory of divisible abelian groups (see e.g.\ \emph{op.\ cit.} D1.1.7). In computer science, (most) relational database schemas can be seen as regular theories over a finite relational signature, as (most) dependencies can be seen as  regular sequents (see e.g.\ \cite{abiteboulhullvianu:95}).  For a regular theory \theory, the category $\modcat{\theory}$ of models and homomorphisms has (small) products and filtered colimits. By the results of Makkai in \cite{makkai:90}, the category \homset{FC}{\modcat{\theory}}{\Sets} of filtered colimit preserving functors is the classifying topos of the theory, that is to say, that
\begin{equation}\label{Equation: Makkai rep}\homset{FC}{\modcat{\theory}}{\Sets}\simeq\classtop\end{equation}
where \classtop\ is the unique, up to equivalence, Grothendieck topos with the property that for any Grothendieck topos \cat{E}, the category of \theory-models in \cat{E} is equivalent to the category of geometric morphisms from \cat{E} to \classtop. (From now on,  ``topos'' will be used to mean Grothendieck topos.) The classifying topos contains a `universal' or `generic' topos model of \theory\ in the form of a regular subcategory, which can be taken to represent \theory\ up to a suitable notion of equivalence. It is shown in \cite{makkai:90} that this subcategory can be characterized up to effective completion as those functors that also preserve products.

This note presents a topological view on Makkai's representation result. In the case of propositional regular theories, which we can think of as meet semi-lattices, it can straightforwardly be seen that a two-valued function of the set of models---i.e.\ filters---of a meet semi-lattice
\[\mathrm{Filt}(S)\to 2\]
is a directed join preserving poset map if and only if it is continuous with respect to the Sierpi\'{n}ski topology on $2$ and the usual `Stone duality' topology on $\mathrm{Filt}(S)$. Thus the free frame $\mathrm{Idl}(S)$ on the meet semi-lattice can simultaneously be represented as the direced colimit preserving poset morphisms to $2$ or as the continuous maps to $2$, and the comparing isomorphism between the latter two is the identify on functions. In the predicate logic case, we know from \cite{butz:98b} that, in addition to (\ref{Equation: Makkai rep}), \classtop\ has a topological representation as equivariant sheaves of a topological groupoid of models and isomorphisms. We show here that the comparison between such a topological representation and Makka's functorial representation is essentially the same as in the propositional case: a set valued functor on the category of models preserves filtered colimits if and only if it is `continuous'---i.e.\ and equivariant sheaf---with respect to a suitable topology on that category. To this end, we cut down to a sufficiently large set of models and equip it with a generalization of the `Stone duality' topology from propositional logic. (The topology is a, in some respects somewhat simpler, variant of the one used in \cite{butz:98b}. It is also used in \cite{APAL:fol}, see Section \ref{Section: Concluding remarks}.) We also equip the set of model-homorphisms with a topology, so that the space of models and the space of homomorphisms form a topological category \topo{H}. We show that the forgetful functor from equivariant sheaves on \topo{H} to set-valued functors $\Sh{\topo{H}}\to [\topo{H},\Sets]$  that forgets the topology factors as an equivalence
\[\Sh{\topo{H}}\simeq \homset{FC}{\topo{H}}{\Sets}\ \to [\topo{H},\Sets]\]
where $\homset{FC}{\topo{H}}{\Sets}\simeq \classtop$ by Makkai's results.
As a consequence we obtain from this comparison a representation of \classtop\ in terms of a topological category of models and homomorphisms, where the restricted topological groupoid of models and isomorphisms classifies a different (non-regular) theory (see Section \ref{Section: Concluding remarks}).

Section \ref{Section:Preliminaries} recalls regular theories and some of their properties and states the representation result displayed in (\ref{Equation: Makkai rep}). This is not explicitly stated in \cite{makkai:90}, but it is a corollary of the results presented there. Section \ref{Section: Sheaves on topological categories of models} presents a topology on the category of models of a regular theory and shows that the category of filtered colimit preserving functors into sets, and hence the classifying topos of the theory, is equivalent to equivariant sheaves on the resulting topological category. Finally, Section \ref{Section: Concluding remarks}  considers decidable coherent and first-order theories.

\section{Preliminaries}
\label{Section:Preliminaries}
\subsection{Regular theories}
\label{Subsection: Regular theories}
Let $\Sigma$ be a (first-order, possibly many-sorted) signature. A formula over $\Sigma$ is \emph{regular} if it is constructed using only the connectives $\top$, $\wedge$, and $\exists$. We consider formulas-in-context $\syntob{\alg{x}}{\phi}$, where the context \alg{x} is a list of distinct variables (implicitly typed) containing (at least) the free variables of $\phi$.  A sequence $\phi\vdash_{\alg{x}}\psi$, where \alg{x} is a context for both $\phi$ and $\psi$, is \emph{regular} if both $\phi$ and $\psi$ are regular formulas. A \emph{regular theory} \theory\ is a deductively closed set of regular sequents. See e.g.\ \cite{elephant1} for further details and a calculus for regular sequents.

Let $\modcat{\theory}$ be the category of \theory-models and homomorphisms for a regular theory \theory. Then $\modcat{\theory}$ has (small) products and filtered colimits, and these can be computed in the category of $\Sigma$-structures and homomorphisms. Furthermore, recall the following notions and properties from \cite{makkai:90}:
\begin{definition}\label{Definition: pp model}
\begin{enumerate}
\item For a \theory-model \alg{M} and regular formula \syntob{\alg{x}}{\phi }, an element $\alg{a}$ is a \emph{generic element} of \alg{M} for \syntob{\alg{x}}{\phi} if $\alg{a}\in\csem{\alg{x}}{\psi}^{\alg{M}}$ implies that the sequent $\phi\vdash_{\alg{x}}\psi$ is provable in \theory.
\item For a tuple $\alg{a}$ of a model \alg{M}, a formula \syntob{\alg{x}}{\phi} is said to \emph{isolate} $\alg{a}$ (in \alg{M}) if $\alg{a}\in\csem{\alg{x}}{\phi }^{\alg{M}}$ and for any model \alg{N} and $\alg{b}\in \csem{\alg{x}}{\phi}^{\alg{N}}$, there exists homomorphism $\alg{f}:\alg{M}\rightarrow\alg{N}$ with $\alg{f}(\alg{a})=\alg{b}$.
\item \alg{M} is a \emph{principal prime} (pp) model if for every tuple $\alg{a}$ of \alg{M} there is a formula isolating $\alg{a}$.
\end{enumerate}
\end{definition}
%
%
\begin{proposition}\label{proposition: universal pp model}
Let a regular formula \syntob{\alg{x}}{\phi} of \theory\ be given. Then there exists a pp model \alg{M} in $M_{\theory}$ and  an element $\alg{a}$ such that \alg{M} has size not larger than $|\Sigma|+\aleph_0$ and \alg{a} is generic for \syntob{\alg{x}}{\phi}. Consequently, \syntob{\alg{x}}{\phi} isolates \alg{a} in \alg{M}.
\begin{proof}The proofs can be found in Section 4  of \cite{makkai:90}.
\end{proof}
\end{proposition}
For a formula \syntob{\alg{x}}{\phi}, refer to a model \alg{M} with element \alg{a} satisfying the properties of Proposition \ref{proposition: universal pp model} as a \emph{universal model} and a \emph{generic element} for \syntob{\alg{x}}{\phi}. (Note that universal models need not be unique up to isomorphism.)

\subsection{The classifying topos of a regular theory}
\label{Subsection: The classifying topos of a regular theory}
Let \theory\ be a regular theory. Let \homset{FC}{\modcat{\theory}}{\Sets} be the category of filtered colimit preserving functors from $\modcat{\theory}$ to $\Sets$ (or, equivalently, of those functors that preserve colimits of directed diagrams, see \cite[1.A.1.5]{adamekandrosicky:94}).  Let $\kappa$ be a regular cardinal strictly larger than $\Sigma$,  let  $\mathrm{Mod}^{\kappa}_{\theory}$ be the essentially small category of strictly smaller than $\kappa$ models, and consider the category $\homset{FC_{\kappa}}{\mathrm{Mod}^{\kappa}_{\theory}}{\Sets}$ of those functors that preserve colimits of diagrams $\cat{D}\rightarrow \mathrm{Mod}^{\kappa}_{\theory}$ such that \cat{D} is a directed poset strictly smaller than $\kappa$. It is straightforward to see that the restriction functor from $\homset{FC}{\modcat{\theory}}{\Sets}$ to $\homset{FC_{\kappa}}{\mathrm{Mod}^{\kappa}_{\theory}}{\Sets}$ is (one half of) an equivalence.
\begin{equation}\label{Equation: Cutting down to kappa models}\homset{FC}{\modcat{\theory}}{\Sets}\simeq \homset{FC_{\kappa}}{\mathrm{Mod}^{\kappa}_{\theory}}{\Sets}\end{equation}
%
%
%

%
Recall (from e.g.\ \cite{elephant1}) that the \emph{syntactic category} \synt{C}{T} of a regular theory is the regular category with the universal property that the category \modin{\theory}{R} of \theory-models  in a regular category \cat{R} is equivalent to the category of regular functors from \synt{C}{T} to \cat{R},
\[\modcat{\theory}{(\cat{R})}\simeq \homset{}{\synt{C}{T}}{\cat{R}}\]
(naturally in \cat{R}). \synt{C}{T} can be constructed from the syntax of \theory\ by taking the objects to be ($\alpha$-equivalence classes) of formulas in context \syntob{\vec{x}}{\phi}, and morphisms  $\syntob{\vec{x}}{\phi}\to\syntob{\vec{y}}{\psi}$ to be (\theory-provable equivalence classes of)  formulas in context \syntob{\vec{x},\vec{y}}{\sigma} such that \theory\ proves that $\sigma$ is a functional relation between $\phi$ and $\psi$. Further, recall (again e.g.\ from \cite{elephant1}) that the \emph{classifying topos} \classtop\ is the topos with the universal property that the category \modin{\theory}{\cat{E}} of \theory-models  in a topos \cat{E} is equivalent to the category of geometric morphisms from \cat{E} to \classtop,
\[\modcat{\theory}{(\cat{E})}\simeq \homset{}{\cat{E}}{\classtop}\]
(naturally in \cat{E}). \classtop\ can be constructed from \synt{C}{T} by taking sheaves on \synt{C}{T} equipped with the regular coverage (see \emph{ibid}.\ for further details). The following emerges  from the results and constructions of \cite{makkai:90}. Since it is not explicitly stated there, a sketch of a proof is included here with references to the relevant results in \cite{makkai:90}.
\begin{theorem}[Makkai]\label{Theorem: Makkais representation theorem}
For \theory\ a regular theory, the category \homset{FC}{\modcat{\theory}}{\Sets} is the classifying topos of \theory,
\[\homset{FC}{\modcat{\theory}}{\Sets}\simeq\classtop\]
\begin{proof} By Proposition 5.2 of \cite{makkai:90}, the \emph{evaluation functor}
\[Ev:\synt{C}{T}\rightarrow \homset{FC}{\modcat{\theory}}{\Sets}\]
which sends an object \syntob{\alg{x}}{\phi} to the `evaluation' or `definable set' functor $\alg{M}\mapsto\csem{\alg{x}}{\phi}^{\alg{M}}$ is regular, full, and faithful. Moreover, it clearly reflects covers, in the sense that if a family of morphisms with common codomain in \synt{C}{T} are sent to a jointly covering family, then there exists a morphism in the family which is a cover in \synt{C}{T}. By Lemma 5.5 of \cite{makkai:90}, the definable set functors (i.e.\ the image of ev) form a generating set for  $\homset{FC}{\modcat{\theory}}{\Sets}$ (note the remark in the proof of 5.3 that only the assumption that the functors preserve filtered colimits is used). The category $\homset{FC}{\modcat{\theory}}{\Sets}$ has small hom-sets, as we see from Equation (\ref{Equation: Cutting down to kappa models}) above.  Finally, since filtered colimits commute with finite limits in \Sets, and limits and colimits are computed `pointwise' in the functor category $[\modcat{\theory},\Sets]$, the category \homset{FC}{\modcat{\theory}}{\Sets} satisfies Giraud's conditions (for a statement of which see e.g.\ \cite{maclane92}), so it is a topos. And by the above, it is equivalent to the category of sheaves on \synt{C}{T} with respect to the regular coverage R, so $\homset{FC}{\modcat{\theory}}{\Sets}\simeq \Sh{\synt{C}{T},R}\simeq\classtop$.
\end{proof}
\end{theorem}

\subsection{Topological categories and sheaves}
\label{Subsection: Topological categories and sheaves}
Consider a topological category, \topo{C}, that is, a category object in the category \alg{Sp} of topological spaces,
\begin{equation}\bfig
%
%
\morphism|m|<750,0>[C_1\times_{C_0} C_1`C_1;\circ]
%
%
\morphism(750,0)|a|/@{>}@<5pt>/<750,0>[C_1`C_0;d]
\morphism(750,0)|m|/@{<-}/<750,0>[C_1`C_0;e]
\morphism(750,0)|b|/@{>}@<-5pt>/<750,0>[C_1`C_0;c]
%
%
\efig\end{equation}
where $\circ$ is composition, $e$ is the identities map, and $d$ and $c$ are the domain and codomain maps, respectively. Recall (see \cite{elephant1}, \cite{moerdijk:88}, \cite{moerdijk:90}) that the objects of the category of \emph{equivariant sheaves}, \Eqsheav{C_1}{C_0}, on \topo{C} are pairs \pair{r:R\rightarrow C_0,\rho} where $r$ is a local homeomorphism---i.e.\ an object of \Sh{C_0}---and $\rho$ is a continuous action, i.e.\ a continuous map
\[\rho:C_{1}\times_{C_{0}}R\to R\]
(with the pullback being along the domain map) such that $r(\rho(f,x))=c(f)$, satisfying a \emph{unit} and a \emph{composition} axiom:
\[ \rho(1_{r(x)},x)=x\ \ \ \ \ \  \rho(g\circ f,x)=\rho(g,\rho(f,x))\]
%
%
A morphism of equivariant sheaves is a morphism of sheaves commuting with the actions. The category, \Eqsheav{C_1}{C_0}, of equivariant sheaves on \topo{C} is a  topos.
%
%
%
A \emph{continuous functor}, or morphism of topological categories, $f:\topo{C}\to\topo{D}$, i.e.\ a morphism of category objects in \alg{Sp}
\[\bfig
\square(0,300)|arra|<1000,300>[C_1\times_{C_0}C_1`D_1\times_{D_0}D_1`C_1`D_1;f_1\times f_1`\circ`\circ`]
\square(0,0)|allb|/>`@{>}@<-8pt>`@{>}@<-8pt>`>/<1000,300>[ C_{1}`D_{1}`C_{0}`D_{0};f_1`d`d`f_0]
\square(0,0)|arrb|/>`@{>}@<8pt>`@{>}@<8pt>`>/<1000,300>[ C_{1}`D_{1}`C_{0}`D_{0};f_1`c`c`f_0]
\square(0,0)|arrb|/>`@{<-}@<-2pt>`@{<-}@<-2pt>`>/<1000,300>[ C_{1}`D_{1}`C_{0}`D_{0};f_1`e`e`f_0]
\efig\]
induces a \emph{geometric morphism}, $f:\Eqsheav{C_1}{C_0}\to\Eqsheav{D_1}{D_0}$, that is, a pair of adjoint functors,
\[\bfig
\morphism|b|/{@{>}@/_1em/}/<1000,0>[\Eqsheav{C_1}{C_0}`\Eqsheav{D_1}{D_0};f_*]
\morphism/{@{<-}@/^1em/}/<1000,0>[\Eqsheav{C_1}{C_0}`\Eqsheav{D_1}{D_0};f^*]
\place(500,0)[\bot]
\efig\]
consisting of a \emph{direct image} functor, $f_*$, and an \emph{inverse  image} functor, $f^*$, where the latter preserves finite limits (and therefore, being a left adjoint, regular logic). A special case is the continuous functor consisting of identity maps to a topological category \topo{C} from the topological category $\topo{C}^d$ consisting of \topo{C} equipped with discrete topologies. In that case the inverse image part of the induced geometric morphism is the forgetful functor
\[u^*:\Eqsheav{C_1}{C_0}\to\Eqsheav{C_1^d}{C_0^d}\simeq\Sets^{\cat{C}}\]
Explicitly, $u^*$ sends an equivariant sheaf, \pair{r:R\rightarrow C_0,\rho}, to the functor which sends an arrow $f:x\rightarrow y$ in $C_1$ to the function $\rho(f,-):r^{-1}(x)\rightarrow r^{-1}(y)$. Note that superscript $d$ will be used also below to indicate discrete or no topology.
%
%
%
%

%
\section{Sheaves on topological categories of models}
\label{Section: Sheaves on topological categories of models}
\subsection{The topological category of countable  models}
\label{Subsection: Topological categories of models}
Fix a (first-order with equality, possibly many-sorted) signature $\Sigma$ and a regular theory \theory\ over $\Sigma$. Fix an infinite set, $\mathds{D}$, of `elements' such that $|\mathds{D}|\geq |\Sigma|$. For simplicity and convenience of presentation, it will be assumed that $\Sigma$ is countable and that $\mathds{D}$ is countably infinite, but the construction and results of this paper does not depend on that assumption. $\mathds{D}$ should be thought of as the `universe' from which we construct $\Sigma$-structures.
Let $\topo{M}_{\theory}^d$ be the category consisting of the set $M^d_{\theory}$ of all $\theory$-models with elements from $\mathds{D}$ (i.e.\ such that underlying sets are contained in $\mathds{D}$), and the set $H^d_{\theory}$ of homomorphisms between them. Clearly, $\topo{M}_{\theory}^d$ is equivalent to the category $\mathrm{Mod}^{\aleph_0}_{\theory}$ of countable \theory-models, so that
\begin{equation}\label{Equation: Equivalences}\homset{FC_{\aleph_0}}{\topo{M}_{\theory}^d}{\Sets}\simeq \homset{FC_{\aleph_0}}{\mathrm{Mod}^{\aleph_0}_{\theory}}{\Sets}\simeq \homset{FC}{\mathrm{Mod}_{\theory}}{\Sets}\simeq \classtop\end{equation}
%
The purpose of this section is to show that the filtered colimit preserving functors on $\topo{M}_{\theory}^d$ are precisely the equivariant sheaves on $\topo{M}_{\theory}^d$ equipped with the following topology.
Certain notational shortcuts will be used: in addition to tuples of variables and elements, boldface will also be used for models and model homomorphisms (as a reminder that if $\Sigma$ is many-sorted, then homomorphisms are families of functions). Typing will usually be left implicit, including in expressions such as $\alg{a}\in\alg{M}$ and $\alg{f}(\alg{a})=\alg{b}$ (for  a model \alg{M} and a homomorphism \alg{f}), trusting that this will cause no confusion.
\begin{definition}\label{Definition> Logical topology}
Let $M_{\theory}$ be the space obtained by equipping $M_{\theory}^d$ with the coarsest topology containing all sets of the following form:
\begin{enumerate}
\item For each sort, $A$, and element, $a\in \mathds{D}$, the set \[\bopen{A,a}:=\cterm{\alg{M}\in M_{\theory}}{a\in \sem{A}^{\alg{M}}}\]
\item For each relation symbol, $R:A_1,\ldots,A_n$, and $n$-tuple, $\alg{a}\in \mathds{D}$, of elements, the set \[\bopen{R,\alg{a}}:=\cterm{\alg{M}\in M_{\theory}}{\alg{a}\in \sem{R}^{\alg{M}}\subseteq \sem{A_1}^{\alg{M}}\times\ldots\times\sem{A_n}^{\alg{M}}}\]
    \textup{(}This extends to nullary relations symbols; if $R$ is a nullary relation symbol, then  $\bopen{R}=\cterm{\alg{M}\in M_{\theory}}{\alg{M}\models R}$.
\item For each function symbol, $f:A_1,\ldots,A_n\rightarrow B$, and n+1-tuple of elements, $\pair{a_1,\ldots,a_n,b}$, the set \[\bopen{f(\alg{a})=b}:=\cterm{\alg{M}\in M_{\theory}}{\sem{f}^{\alg{M}}(\alg{a})=b}\] \textup{(}This extends in the obvious way to include nullary function symbols, i.e.\ constants\textup{)}.
\end{enumerate}
Let $H_{\theory}$ be the space obtained by equipping $H_{\theory}^d$  with the coarsest topology such that the \emph{domain} and \emph{codomain} functions $d,c: H_{\theory}\rightrightarrows M_{\theory}$ are both continuous and such that all sets of the following form are open:
\begin{enumerate}[(i)]
\item For each sort, $A$, and pair of elements, $a,b\in \mathds{D}$, the set
 \[\bopen{A,a\mapsto b}:=\cterm{\alg{f}\in H_{\theory}}{a\in \sem{A}^{d(\alg{f})}\ \textnormal{and}\ f_{A}(a)=b}\]
\end{enumerate}
\end{definition}
%
%
\begin{lemma}\label{Lemma: C continuous category}
The spaces $M_{\theory}$ and $H_{\theory}$ form a topological category $\topo{M}_{\theory}$, i.e.\ a category object in the category, \alg{Sp}, of topological spaces and continuous functions:
\[\bfig
%
%
\morphism<750,0>[H_{\theory}\times_{M_{\theory}} H_{\theory}`H_{\theory};\circ]
%
%
\morphism(750,0)|a|/@{>}@<5pt>/<750,0>[H_{\theory}`M_{\theory};d]
\morphism(750,0)|m|/@{<-}/<750,0>[H_{\theory}`M_{\theory};e]
\morphism(750,0)|b|/@{>}@<-5pt>/<750,0>[H_{\theory}`M_{\theory};c]
\efig\]
\begin{proof} Straightforward.
\end{proof}
\end{lemma}
Call the topologies of Definition \ref{Definition> Logical topology} the \emph{logical} topologies. Note that they are defined in terms of the structure $\Sigma$ and not the theory \theory\ (or the objects in \synt{C}{T}). However, it is convenient to note that a basis for the topology of $M_{\theory}$ can be given in terms of the formulas over $\Sigma$.
\begin{lemma}\label{Lemma: The logical topologies are all the same}
Sets of the form
\[\bopen{\syntob{\alg{x}}{\phi},\alg{a}}=\cterm{\alg{M}\in M_{\theory}}{\alg{a}\in \csem{\alg{x}}{\phi}^{\alg{M}}}\]
form basis for the logical topology on $M_{\theory}$, with $\phi$ ranging over all regular formulas over $\Sigma$.
\begin{proof} By a straightforward induction.
\end{proof}
\end{lemma}
Accordingly, a basic open set of $H_{\theory}$ can be written, displaying a `domain', `preservation', and `codomain' condition, as
\[\left(\begin{array}{c}
\syntob{\alg{x}}{\phi},\alg{a}  \\
\alg{z}:\alg{b} \mapsto \alg{c}   \\
 \syntob{\alg{y}}{\psi},\alg{d}
\end{array}\right)\]
where the sorts for the preservation condition are given by the implicitly typed variables \alg{z}.
\begin{remark}\label{Remark: Decriptive set theory}
The topology on countable models for a regular theory defined in \ref{Definition> Logical topology} may remind the reader of certain well-studied topologies used to form Polish spaces of models on a fixed countable domain for first-order theories, see e.g.\ \cite{grzemostryll:61} or the expository \cite{hjort:groupactionsandmodels}. Now, the space of countable models of \ref{Definition> Logical topology} is not Polish, but by considering (presentations of) the basic open sets in a completely prime filter of open sets, one can straightforwardly see that it is sober. See Section \ref{Section: Concluding remarks} for more on \ref{Definition> Logical topology} and first-order theories.
\end{remark}
%
%
%
%
%
%
%
%
%
%
%
\subsection{Equivariant sheaves on models}
\label{Subsection: Equivariant sheaves on MT}
Consider the topos \Eqsheav{H_{\theory}}{M_{\theory}} and the forgetful functor
\[u^*:\Eqsheav{H_{\theory}}{M_{\theory}}\rightarrow \Eqsheav{H^d_{\theory}}{M^d_{\theory}}\simeq\Sets^{\topo{M}_{\theory}^d}\]
(using superscript $d$ to indicate discrete). Recall from Section \ref{Section:Preliminaries} the `evaluation' functor $Ev:\synt{C}{T}\rightarrow \Sets^{\topo{M}^d}$ which sends an object \syntob{\alg{x}}{\phi} to the functor defined by $\alg{M}\mapsto\csem{\alg{x}}{\phi}^{\alg{M}}$ (extend to morphisms in \synt{C}{T} and homomorphisms of models in the expected way). Call the objects in the image of $Ev$ \emph{definable set functors}. By the results of \cite{makkai:90}, the evaluation functor is regular, full and faithful, and factors through the subcategory of filtered colimiting preserving functors
\[Ev:\synt{C}{T}\to\homset{FC_{\aleph_0}}{\topo{M}_{\theory}^d}{\Sets}\embedd \Sets^{\topo{M}^d}\]
It factors through  \Eqsheav{H_{\theory}}{M_{\theory}} in the following way. First, for an object $\syntob{\alg{x}}{\phi}$ of \synt{C}{T}, represent the corresponding definable set functor in \Eqsheav{H^d_{\theory}}{M^d_{\theory}} as the set
\[ \sox{\alg{x}}{\phi}=\cterm{\pair{\alg{M},\alg{a}}}{\alg{M}\in M_{\theory},\alg{a}\in \csem{\alg{x}}{\phi}^{\alg{M}}}\to^{\pi} M_{\theory} \]
over $M_{\theory}$, where the $\pi$ projects out the model, together with the action
\[\theta_{\syntob{\alg{x}}{\phi}}(\pair{\alg{M},\alg{a}},\alg{f}:\alg{M}\rightarrow \alg{N})=\pair{\alg{N},\alg{f}(\alg{a})}\]
The subscript on $\theta$ will usually be left implicit.
\begin{definition}\label{Definition: Logical topology on sheaves}
For an object $\syntob{\alg{x}}{\phi}$ of \synt{C}{T}, the \emph{logical topology} on the set \sox{\alg{x}}{\phi} is the coarsest such that the projection $ \csem{\alg{x}}{\phi}_{M_{\theory}}\rightarrow M_{\theory}$ is continuous and such that for every list of elements $\alg{a}\in\mathds{D}$ of the same length as \alg{x}, the image of the map $\bopen{\syntob{\alg{x}}{\phi},\alg{a}}\rightarrow \sox{\alg{x}}{\phi}$ defined by $\alg{M}\mapsto \pair{\alg{M},\alg{a}}$ is open.
\end{definition}
\begin{lemma}\label{lemma: Logical topology on sheaves}
For an object $\syntob{\alg{x}}{\phi}$ of \synt{C}{T}, a basis for the \emph{logical topology} on the set $\sox{\alg{x}}{\phi}$ is given by sets of the form
\[
\bopen{\syntob{\alg{x},\alg{y}}{\psi}, \alg{b}}:= \cterm{\pair{\alg{M},\alg{a}}}{ \alg{a},\alg{b}\in
\csem{\alg{x},\alg{y}}{\phi\wedge\psi}^{\alg{M}}}
\]
(where $\alg{b}$ is of the same length as $\alg{y}$)
\begin{proof}Note that the (open) image of the map $\bopen{\syntob{\alg{x}}{\phi},\alg{a}}\rightarrow \sox{\alg{x}}{\phi}$ defined by $\alg{M}\mapsto \pair{\alg{M},\alg{a}}$ can be written as \bopen{\syntob{\alg{x},\alg{y}}{ \alg{x}=\alg{y}}, \alg{a}}. In general,  \bopen{\syntob{\alg{x},\alg{y}}{\psi}, \alg{b}} is an open set: for if $\pair{\alg{M},\alg{a}}\in\bopen{\syntob{\alg{x},\alg{y}}{\psi}, \alg{b}}$, then $\pair{\alg{M},\alg{a}}\in \bopen{\syntob{\alg{x},\alg{y}}{ \alg{x}=\alg{y}}, \alg{a}}\cap \pi^{-1}(\bopen{\syntob{\alg{x},\alg{y}}{\psi}, \alg{a},\alg{b}})\subseteq \bopen{\syntob{\alg{x},\alg{y}}{\psi}, \alg{b}}$. It is clear that such sets form a basis.
\end{proof}
\end{lemma}
\begin{proposition}\label{Proposition: DC, Md factors as M}
The mapping $\syntob{\alg{x}}{\phi}\mapsto \pair{\sox{\alg{x}}{\phi}, \theta}$ defines the object part of a functor $\cat{M}:\synt{C}{T}\to\Eqsheav{H_{\theory}}{M_{\theory}}$
which factors the evaluation functor through the forgetful inverse image functor (up to natural isomorphism),
\[\bfig
\Vtriangle/>`<-`<-/<350,350>[\Eqsheav{H_{\theory}}{M_{\theory}}`\Sets^{\topo{M}_{\theory}^d}`\synt{C}{T};u^*`\cat{M}`\mathrm{ev}]
\efig\]
Consequently, $\cat{M}$ is regular, full, faithful, and cover reflecting (with respect to the regular coverage on \synt{C}{T}).
\begin{proof}
It is clear from Definition \ref{Definition: Logical topology on sheaves} and Lemma \ref{lemma: Logical topology on sheaves} that the projection $\pi:\sox{\alg{x}}{\phi}\rightarrow X_{\theory}$ is a local homeomorphism.
Also, given an arrow
\[\syntob{\alg{x},\alg{y}}{\sigma}:\syntob{\alg{x}}{\phi}\to
\syntob{\alg{y}}{\psi}\]
in \synt{C}{T}, the function $f_{\sigma}=\cat{M}(\sigma): \sox{\alg{x}}{\phi}\rightarrow \sox{\alg{y}}{\psi}$ is
continuous. For given a basic open
$\bopen{\syntob{\alg{y},\alg{z}}{\xi}, \alg{c}} \subseteq \sox{\alg{y}}{\psi}$, then
\[
f_{\sigma}^{-1}\left(\bopen{\syntob{\alg{y},\alg{z}}{\xi}, \alg{c}} \right) =
\bopen{\syntob{\alg{x},\alg{z}}{\fins{\alg{y}} \sigma \wedge \xi}, \alg{c}}
\]
Next, the action $\theta_{\syntob{\alg{x}}{\phi}}$ is continuous: let a basic open
$U= \bopen{\syntob{\alg{x},\alg{y}}{\psi}, \alg{b}}\subseteq \sox{\alg{x}}{\phi}$ be
given, and suppose $\theta(\alg{f},\pair{\alg{M},\alg{a}})= \pair{\alg{N},\alg{c}}\in U$ for
$\alg{M},\alg{N}\in M_{\theory}$ and $f:\alg{M}\rightarrow \alg{N}$ in $H_{\theory}$. Then we can specify an
open neighborhood around \pair{\alg{f},\pair{\alg{M},\alg{a}}} which $\theta$ maps into $U$ as:
\[
\pair{\alg{f},\pair{\alg{M},\alg{a}}}\in \left(\begin{array}{c}
-  \\
  \alg{x}:\alg{a} \mapsto \alg{c}   \\
 \bopen{\syntob{\alg{x},\alg{y}}{\psi}, \alg{c},\alg{b}}
\end{array}\right)
 \times_{M_{\theory}}
\bopen{\syntob{\alg{x},\alg{z}}{\alg{x}=\alg{z}}, \alg{a}}
\]
It is clear that \cat{M} factors the evaluation functor through the forgetful functor, and---since the latter is a conservative inverse image functor---that the fact that the evaluation functor is regular, full, faithful, and cover reflecting implies that \cat{M} has the same properties.
\end{proof}
\end{proposition}
The key observation is now that $R$ is an equivariant sheaf over $\topo{M}_{\theory}$ then the underlying functor $u^*(R)$ must preserve filtered colimits. We need the following technical lemma, stating that for an element $x$ in an equivariant sheaf over $\topo{M}_{\theory}$ there exists a section around $x$ which is, in a certain sense, `well-behaved'.
\begin{lemma}\label{Lemma: Well behaved section}
For any object
\[\left\langle\bfig\morphism<300,0>[R`M_{\theory};r]\efig,\rho\right\rangle\]
in \Eqsheav{H_{\theory}}{M_{\theory}},
%
%
and any element $x\in R$, there exists a basic open $\bopen{\syntob{\alg{x}}{\phi},\alg{a}}\subseteq
M_{\theory}$ and a section $v: \bopen{\syntob{\alg{x}}{\phi},\alg{a}}\rightarrow R$ containing $x$ such
that for any $\alg{f}:\alg{M}\rightarrow \alg{N}$ in $H_{\theory}$ such that $\alg{M}\in
\bopen{\syntob{\alg{x}}{\phi},\alg{a}}$ and $\alg{f}(\alg{a})=\alg{a}$ (so that $\alg{N}$ is also in
$\bopen{\syntob{\alg{x}}{\phi},\alg{a}}$), we have $\rho(\alg{f},v(\alg{M}))= v(\alg{N})$.
\begin{proof}
Given $x\in R$, choose a section $s:\bopen{\syntob{\alg{y_1}}{\psi},\alg{b_1}}\rightarrow R$ such that $x\in
s(\bopen{\syntob{\alg{y_1}}{\psi},\alg{b_1}})$. Pull the open set
$s(\bopen{\syntob{\alg{y_1}}{\psi},\alg{b_1}})$ back along the continuous action $\rho$,
\[
\bfig \square|allb|/->` >->`
>->`->/<600,300>[V`s(\bopen{\syntob{\alg{y_1}}{\psi},\alg{b_1}})`H_{\theory}\times_{M_{\theory}} R`R;`\subseteq`
\subseteq`\rho ] \place(75,250)[\spbangle]  \efig
\]
to obtain an open set $V$ containing \pair{1_{r(x)},x}. Since $V$ is open, we can find a box $W$ of basic opens
around \pair{1_{r(x)},x} contained in $V$:
\[
\pair{1_{r(x)},x}\in W= \left(\begin{array}{c}
\syntob{\alg{y_2}}{\xi},\alg{b_2}  \\
  \alg{z}:\alg{k} \mapsto \alg{k}   \\
 \syntob{\alg{y_3}}{\eta},\alg{b_3}
\end{array}\right)
\times_{M_{\theory}} v'(\bopen{\syntob{\alg{y_4}}{\theta},\alg{b_4}})\subseteq V
\]
where $v'$ is a section $v':\bopen{\syntob{\alg{y_4}}{\theta},\alg{b_4}}\rightarrow R$ with $x$ in its
image. Notice that the preservation condition of $W$,  i.e.\ ($\alg{z}:\alg{k} \mapsto \alg{k}  $), must have the
same elements on both the domain and the codomain side, since it is satisfied by $1_{r(x)}$. Now, restrict $v'$ to the
subset
\[
U:=\bopen{\syntob{\alg{y_2}, \alg{z}, \alg{y_3}, \alg{y_4}}{\xi \wedge \eta
\wedge \theta},\alg{b_2}, \alg{k}, \alg{b_3}, \alg{b_4}}
\]
to obtain a section $v=v'\upharpoonright_{U}:U\rightarrow R$. Notice that $x\in v(U)$. Furthermore,
$v(U)\subseteq s(\bopen{\syntob{\alg{y_1}}{\psi},\alg{b_1}})$, for if $v(\alg{M})\in v(U)$, then
$\pair{1_{\alg{M}},v(\alg{M})}\in W$, and so $\rho(\pair{1_{\alg{M}},v(\alg{M})})=v(\alg{M})\in
s(\bopen{\syntob{\alg{y_1}}{\psi},\alg{b_1}})$. Finally, if $\alg{M}\in U$ and $\alg{f}:\alg{M}\rightarrow
\alg{N}$ is a morphism in $H_{\theory}$ such that
\[ \alg{f}(\alg{b_2}, \alg{k}, \alg{b_3}, \alg{b_4}) =\alg{b_2}, \alg{k}, \alg{b_3}, \alg{b_4}\]
then $\pair{\alg{f},v(\alg{M})} \in W$, and so $\rho(\alg{f},v(\alg{M}))\in
s(\bopen{\syntob{\alg{y_1}}{\psi},\alg{b_1}})$. But we also have $v(\alg{N})\in v(U) \subseteq
s(\bopen{\syntob{\alg{y_1}}{\psi},\alg{b_1}})$, and $r(\rho(\alg{f},v(\alg{M})))=r(v(\alg{N})$,  so
$\rho(\alg{f},v(\alg{M}))=v(\alg{N})$. Hence $v:\bopen{\syntob{\alg{y_2}, \alg{z}, \alg{y_3}, \alg{y_4}}{\xi \wedge \eta
\wedge \theta},\alg{b_2}, \alg{k}, \alg{b_3}, \alg{b_4}}\to R$ is the sought for section with respect to the given $x\in R$.
\end{proof}
\end{lemma}
\begin{lemma}\label{Lemma: Eq sheaves preserve filtered colimits}
The forgetful functor $u^*:\Eqsheav{H_{\theory}}{M_{\theory}}\rightarrow \Sets^{\topo{M}_{\theory}^d}$ factors through $\homset{FC}{\topo{M}_{\theory}^d}{\Sets}\hookrightarrow\Sets^{\topo{M}_{\theory}^d}$.
\begin{proof}
For convenience of presentation, we think of $\Sigma$ as single-sorted.
We must show that the underlying functor of an equivariant sheaf \pair{r:R\rightarrow M_{\theory},\rho} preserves colimits of countable, directed diagrams. Let
$\alg{M}_{-}:\cat{D}\rightarrow \topo{M}_{\theory}$ be such a diagram, with \alg{M} the colimit,  $\alg{f}_d:\alg{M}_d\rightarrow \alg{M}$ the morphisms of the colimiting cone, and $\alg{g}_{d,d'}:\alg{M}_d\rightarrow\alg{M}_{d'}$ the morphisms of the diagram.
\[\bfig
\Atriangle/<-`<-`/<500,550>[\alg{M}\cong\colimit_{d\in \cat{D}}\alg{M}_d`\alg{M}_a`\alg{M}_b;\alg{f}_{a}`\alg{f}_{b}`]
\morphism(500,200)<0,350>[\alg{M}_c`\alg{M}\cong\colimit_{d\in \cat{D}}\alg{M}_d;\alg{f}_c]
\Atriangle|bbb|/<-`<-`/<500,200>[\alg{M}_c`\alg{M}_a`\alg{M}_b;\alg{g}_{a,c}`\alg{g}_{b,c}`]
\efig\]
We lose no generality in making the following assumption: for all elements $a\in\alg{M}$ there exists $d\in\cat{D}$ such that $a\in\alg{M}_d$; and for all $a\in\alg{M}$ and all $d\in\cat{D}$, if $a\in \alg{M}_d$ then $\alg{f}_d(a)=a$ and for all $d'\geq d$,  $\alg{g}_{d,d'}(a)=a$.

\noindent First, the net (or generalized sequence) $\alg{M}_{-}:\cat{D}\rightarrow M_{\theory}$ converges to \alg{M}: For if $\alg{M}\in \bopen{\syntob{\alg{x}}{\phi},\alg{a}}$ then, by the construction of directed colimits of structures and the assumption above, there exists $d\in \cat{D}$ such that for all $d'\geq d$ we have   $\alg{M}_{d'}\in \bopen{\syntob{\alg{x}}{\phi},\alg{a}}$. Similarly, for all $d\in \cat{D}$, the net $\alg{g}_{d,-}:\uparrow(d)\rightarrow H_{\theory}$ defined by $d'\mapsto \alg{g}_{d,d'}$ converges to $\alg{f}_d$.

\noindent Second, $\cterm{\rho(\alg{f}_{d},-):r^{-1}(\alg{M}_{d})\rightarrow r^{-1}(\alg{M})}{d\in \cat{D}}$
is a jointly surjective family of functions: For, given $x\in r^{-1}(\alg{M})$, we can chose, by Lemma \ref{Lemma: Well behaved section}, a section $v:\bopen{\syntob{\alg{y}}{\phi},\alg{a}}\rightarrow R$ with $x\in v(\bopen{\syntob{\alg{y}}{\phi},\alg{a}})$ such that for any $\alg{K}\in \bopen{\syntob{\alg{y}}{\phi},\alg{a}}$ and $\alg{g}:\alg{K}\rightarrow\alg{L}$ such that $\alg{g}(\alg{a})=\alg{a}$, we have that $\rho(\alg{g},v(\alg{K}))=v(\alg{L})$. Since the net $\cterm{\alg{M}_d}{d\in \cat{D}}$ converges to \alg{M}, we can find $d\in \cat{D}$ such that for all $d'\geq d$ we have that $\alg{M}_{d'}\in \bopen{\syntob{\alg{y}}{\phi},\alg{a}}$. But then $\alg{f}_{d'}(\alg{a})=\alg{a}$, so that $\rho(\alg{f}_{d'},v(\alg{M}_{d'}))=v(\alg{M})=x$.

\noindent Finally, for given $d\in\cat{D}$, if $x,y\in r^{-1}(\alg{M}_{d})$  and $\rho(\alg{f}_{d},x)=\rho(\alg{f}_{d},y)$, then there exists $d'\geq d$ such that $\rho(\alg{g}_{d,d'},x)=\rho(\alg{g}_{d,d'},y)$: Since the net $\cterm{ \pair{\alg{g}_{d,d'},x}}{d'\geq d}$ converges to $\pair{\alg{f}_{d},x}$ in $H_{\theory}\times_{M_{\theory}}R$ and the action $\rho:H_{\theory}\times_{M_{\theory}}R\rightarrow R$ is continuous, the net $\cterm{ \rho(\alg{g}_{d,d'},x)}{d'\geq d}$ converges to $\rho(\alg{f}_{d},x)$ in $R$, and likewise $\cterm{ \rho(\alg{g}_{d,d'},y)}{d'\geq d}$ converges to $\rho(\alg{f}_{d},y)$ . But $\rho(\alg{f}_{d},x)=\rho(\alg{f}_{d},y)$, so for any section $w:U\rightarrow R$ containing $\rho(\alg{f}_{d},x)$ there must, then, exist $d'\geq d$ such that $\rho(\alg{g}_{d,d'},x),\rho(\alg{g}_{d,d'},y)\in w(U)$, which means that $\rho(\alg{g}_{d,d'},x)=w(\alg{M}_{d'})=\rho(\alg{g}_{d,d'},y)$.
We conclude that
\[r^{-1}(\alg{M})= \colimit_{d\in \cat{D}}r^{-1}(\alg{M}_d)\] in \Sets.
\end{proof}
\end{lemma}

We claim that the forgetful functor $u^*:\Eqsheav{H_{\theory}}{M_{\theory}}\rightarrow \homset{FC}{\topo{M}_{\theory}^d}{\Sets}$ is one half of an equivalence of categories. Clearly, $u^*$ is faithful and, since the definable set functors generate \homset{FC}{\topo{M}_{\theory}^d}{\Sets}, essentially surjective. We conclude the argument by showing that the definables also form a generating set in \Eqsheav{H_{\theory}}{M_{\theory}}, helping ourselves to another of Makkai's results along the way.

\begin{definition}Given a basic open $\bopen{\syntob{\alg{x}}{\phi},\alg{a}}\subseteq M_{\theory}$ we say that it is (presented) in \emph{reduced form} if $\alg{a}$ has the property with respect to the type of $\alg{x}$ that $a_i=a_j$ implies that either $i=j$ or the sort of $x_i$ is distinct from the sort of $x_j$.
\end{definition}
Clearly, given a presentation of a basic open set we can always produce, by a straightforward substitution, a presentation of the same set which is in reduced form.
\begin{lemma}\label{Lemma: Basic opens of X generate representables}
Let $\bopen{\syntob{\alg{x}}{\phi},\alg{a}}$ be a non-empty basic open of $M_{\theory}$ in reduced form. Then the sheaf \pair{\sox{\alg{x}}{\phi}\rightarrow M_{\theory},\theta} is the closure under $\theta$ of the image of the continuous section
\[v:\bopen{\syntob{\alg{x}}{\phi},\alg{a}} \to \sox{\alg{x}}{\phi}\]
defined by $\alg{M}\mapsto\pair{\alg{M},\alg{a}}$.
\begin{proof}
The displayed section is continuous by Definition \ref{Definition: Logical topology on sheaves}. By Lemma \ref{proposition: universal pp model} there exists a universal model \alg{M} for \syntob{\alg{x}}{\phi}, and we can assume that \alg{a} is a generic element.
\end{proof}
\end{lemma}
Moreover, elements of sheaves `have support' in the following sense of Makkai's.
\begin{lemma}[Makkai]\label{Lemma: Makkai support}
Let $F$ be an object in \homset{FC}{\topo{M}_{\theory}^d}{\Sets}, and let $y\in F(\alg{M})$. Then there exists \syntob{\alg{x}}{\phi} and a list of elements $\alg{a}\in\mathds{D}$ of the same length as \alg{x} such that $\alg{a}\in\csem{\alg{x}}{\phi}^{\alg{M}}$ and for all $\alg{h}_1,\alg{h}_2:\alg{M}\rightrightarrows \alg{N}$ in $\topo{M}_{\theory}$, if $\alg{h}_1(\alg{a})=\alg{h}_2(\alg{a})$ then $F(\alg{h}_1)(y)=F(\alg{h}_2)(y)$.
\begin{proof}See the proof of Lemma 5.3 in \cite{makkai:90}.
\end{proof}
\end{lemma}
%
%
%
%
\begin{lemma}\label{Lemma: Definables are generators} The definable objects form a generating set for $\Eqsheav{H_{\theory}}{M_{\theory}}$.
\begin{proof}
First, let \alg{M} be a pp model, \pair{r:R\rightarrow M_{\theory},\rho} an object in \Eqsheav{H_{\theory}}{M_{\theory}}, and $x\in R$ an element such that $r(x)=\alg{M}$. We show that there exists a morphism with definable domain \[\tilde{\tau}:\cat{M}(\syntob{\alg{x}}{\phi})\to \pair{r:R\rightarrow M_{\theory},\rho}\] such that $x$ is in the image of $\tilde{\tau}$.
We largely follow the construction of \cite[Lemma 5.4]{makkai:90}. Choose a section $\tau:U\rightarrow R$ containing $x$ according to Lemma \ref{Lemma: Well behaved section} and a basic open set $V$ containing \alg{M} according to Lemma \ref{Lemma: Makkai support}, and restrict $\tau$ to $U\cap V$ to obtain a section of the form
\[\tau:\bopen{\syntob{\alg{x}}{\phi}, \alg{a}}\to R\]
containing $x$. We may suppose that  \bopen{\syntob{\alg{x}}{\phi},\alg{a}} is on reduced form. Since \alg{M} is pp, we may also assume that \syntob{\alg{x}}{\phi} isolates $\alg{a}$ (since if $\psi$ is the isolating formula, we can restrict $\tau$ to $\phi\wedge\psi$). Then for any $\alg{N}\in\bopen{\syntob{\alg{x}}{\phi}, \alg{a}}$, there is a homomorphism $\alg{f}:\alg{M}\rightarrow\alg{N}$ such that $\alg{f}(\alg{a})=\alg{a}$, and by Lemma \ref{Lemma: Well behaved section}, we have $\rho(\alg{f},x)=\tau(\alg{N})$. Now, by Lemma \ref{Lemma: Basic opens of X generate representables}, we have a section
\[v:\bopen{\syntob{\alg{x}}{\phi}, \alg{a}}\to \sox{\alg{x}}{\phi}\]
sending $\alg{N}\in\bopen{\syntob{\alg{x}}{\phi}, \alg{a}}$ to \pair{\alg{N},\alg{a}}, such that \sox{\alg{x}}{\phi} is the closure under the action of $v(\bopen{\syntob{\alg{x}}{\phi}, \alg{a}})$. We can therefore define a function
\[\tilde{\tau}: \sox{\alg{x}}{\phi}\to R\]
by choosing, for $\pair{\alg{N},\alg{b}}\in \sox{\alg{x}}{\phi}$ a homomorphism $\alg{f}:\alg{M}\rightarrow \alg{N}$ such that $\alg{f}(\alg{a})=\alg{b}$ and setting $\tilde{\tau}(\pair{\alg{N},\alg{b}}=\rho(\alg{f},x)$. Then $\tilde{\tau}$ is well defined by Lemma \ref{Lemma: Makkai support}, and it clearly commutes with the actions.

\noindent Remains to show that $\tilde{\tau}:\sox{\alg{x}}{\phi}\to R$ is continuous. First, we show that
\[\theta:H_{\theory}\times_{M_{\theory}}\sox{\alg{x}}{\phi}\to \sox{\alg{x}}{\phi}\]
is open; let a non-empty basic open
\[
U=
\left(\begin{array}{c}
\syntob{\alg{y_1}}{\psi},\alg{b_1}  \\
  \alg{y_2}:\alg{b_2} \mapsto \alg{c}   \\
 \syntob{\alg{y_3}}{\vartheta},\alg{b_3}
\end{array}\right)
\times_{M_{\theory}} \bopen{\syntob{\alg{x},\alg{y_4}}{\xi},\alg{b_4} }\subseteq H_{\theory}\times_{M_{\theory}}\sox{\alg{x}}{\phi}
\]
be given. Consider the list $\alg{b_1},\alg{b_2},\alg{b_4}$. Write \syntob{\alg{y_1},\alg{y_2},\alg{y_4}}{\chi} for the formula that expresses the expressible identities that occur in this list, i.e.\ the conjunction of formulas $y_{i_k}=y_{j_l}$ for each $b_{i_k}=b_{j_l}$ such that the type of $y_{i_k}$ is the same as the type of $y_{j_l}$. Now, consider
\[V:=\bopen{\syntob{\alg{x},\alg{y_2},\alg{y_3}}{\fins{\alg{y_1},\alg{y_4}}\chi\wedge\psi\wedge \vartheta\wedge \xi},\alg{c}*\alg{b_3}}\subseteq \sox{\alg{x}}{\phi}\]
Clearly, $\theta(U)\subseteq V$. Conversely, if  $\pair{\alg{N},\alg{m}}\in V$, then we can consider a suitable universal model and generic element for \syntob{\alg{x},\alg{y}_2}{\fins{\alg{y}_1,\alg{y}_4}\chi\wedge\psi\wedge \xi}
to find an member of $U$ which is sent to \pair{\alg{N},\alg{m}}, so $V\subseteq \theta(U)$.

\noindent Second, we  show that $\tilde{\tau}:\sox{\alg{x}}{\phi}\to R$ is continuous. Note that by the construction of $\tilde{\tau}$, the  triangle,
\[\bfig
\Vtriangle/>`<-`<-/[\sox{\alg{x}}{\phi}`R`\bopen{\syntob{\alg{x}}{\phi},\alg{a}};\tilde{\tau}`v`\tau]
\efig\]
commutes, and so the restriction of $\tilde{\tau}$,to the open set $v(\bopen{\syntob{\alg{x}}{\phi},\alg{a}})\subseteq \sox{\alg{x}}{\phi}$ is a homeomorphism
\[v(\bopen{\syntob{\alg{x}}{\phi},\alg{a}})\to^{\cong}{\tau}(\bopen{\syntob{\alg{x}}{\phi},\alg{a}})\]
Consider the following commuting diagram,
\[\bfig
\square(1200,0)<1400,500>[H_{\theory}\times_{M_{\theory}}\sox{\alg{x}}{\phi}`\sox{\alg{x}}{\phi}`H_{\theory}\times_{M_{\theory}}R`R;\theta`1_{H_{\theory}}\times\tilde{\tau}`\tilde{\tau}`\rho]
\square<1200,500>[H_{\theory}\times_{M_{\theory}}v(\bopen{\syntob{\alg{x}}{\phi},\alg{a}})`H_{\theory}\times_{M_{\theory}}\sox{\alg{x}}{\phi}` H_{\theory}\times_{M_{\theory}}{\tau}(\bopen{\syntob{\alg{x}}{\phi},\alg{a}})`H_{\theory}\times_{M_{\theory}}R; \subseteq`\cong``\subseteq]
\efig\]
and note that the top composite morphism is a surjection. Since it is also open, $\tilde{\tau}$ is continuous.

Now for the general case. Let $\pair{r:R\rightarrow M_{\theory},\rho}$ in \Eqsheav{H_{\theory}}{M_{\theory}} and $x\in R$ be given, with $r(x)=\alg{M}$. By \cite[Proposition 4.4]{makkai:90}, we can write \alg{M} as a colimit of a directed diagram of pp models in $\topo{M}_{\theory}$. By Lemma \ref{Lemma: Eq sheaves preserve filtered colimits} $u^*(\pair{r:R\rightarrow M_{\theory},\rho})$ preserves directed colimits, so we can choose a pp model and a homomorphism $\alg{f}_i:\alg{M}_i\rightarrow\alg{M}$ and an element $y\in r^{-1}(\alg{M}_i)$ such that $\rho(\alg{f}_i,y)=x$. Then there exists a definable object and a morphism $\tilde{\tau}:\sox{\alg{x}}{\phi}\rightarrow R$ such that $y$ is in the image of $\tilde{\tau}$. But then so is $x$.
\end{proof}
\end{lemma}
\begin{theorem}\label{Theorem: Representation result}
The forgetful functor $u^*:\Eqsheav{H_{\theory}}{M_{\theory}}\rightarrow \homset{FC}{\topo{M}_{\theory}^d}{\Sets}$ is one half of an equivalence of categories,
\[\Eqsheav{H_{\theory}}{M_{\theory}}\simeq \homset{FC}{\topo{M}_{\theory}^d}{\Sets} \simeq \classtop\]
\begin{proof} By Proposition \ref{Proposition: DC, Md factors as M} the functor $\cat{M}:\synt{C}{T}\rightarrow \Eqsheav{H_{\theory}}{M_{\theory}}$ is regular, full, faithful and cover reflecting w.r.t.\ the regular coverage on \synt{C}{T}. By Lemma \ref{Lemma: Definables are generators}, the image of \cat{M} is generating set. It follows that \Eqsheav{H_{\theory}}{M_{\theory}} is equivalent to the topos of sheaves on \synt{C}{T} equipped with the regular coverage, and therefore equivalent to \classtop, and therefore to \homset{FC}{\topo{M}_{\theory}^d}{\Sets}. Since the forgetful functor is a inverse image functor which restricts to an equivalence of universal \theory-models, it is one half of that equivalence.
\end{proof}
\end{theorem}

\section{Coherent theories}
\label{Section: Concluding remarks}

Makkai's semantic representation of the classifying topos of a regular theory as the filtered colimit preserving functors on models is an elegant and simple construction which only uses structure which is intrinsic to the category of models. Moreover, it characterizes the effective completion of the syntactic category \synt{C}{T} in this topos as the product preserving functors, supplementing the characterization of that category as the stably supercompact objects (see e.g.\ \cite{elephant1}). In this respect, the characterization of the classifying topos as equivariant sheaves on the category of models equipped with (extrinsic) topological structure holds no advantage. Passing to coherent theories, however, there is in general no structure intrinsic to the category of models of a theory in terms of which the classifying topos can be constructed. Makkai's solution is to equip this category with extra structure in terms of ultra-products (see \cite{makkai:87}, \cite{makkai93}). The topological approach presented here, however, extends without change to that case; it is shown in \cite{preprint:fol} (building on results in \cite{butz:98b}) that the classifying topos for a coherent theory with inequality predicates (such as a Morleyized classical first-order theory) can be represented as equivariant sheaves on the  groupoid of models and isomorphisms equipped with the topologies of Definition \ref{Definition> Logical topology}. As a consequence, the classifying topos of such a theory \theory\ can be represented as equivariant sheaves on the topological category $\topo{H}_{\theory}$ defined as in Section \ref{Subsection: Topological categories of models}. We state this and sketch the proof.
\begin{theorem}\label{Theorem: Representation of coherent theories}
Let \theory\ be a coherent theory such that, for all homomorphisms $\alg{h}:\alg{M}\rightarrow \alg{N}$ between \theory-models, the underlying component functions of \alg{h} are 1--1. Let  $\topo{H}_{\theory}$ be the topological category of \theory-models and homomorphisms. Then
\[\classtop\simeq\Eqsheav{H_{\theory}}{M_{\theory}}.\]
\begin{proof} Let $I_{\theory}\subseteq H_{\theory}$ be the subspace of isomorphisms, and $\topo{I}_{\theory}$ the resulting topological groupoid of \theory-models and isomorphisms. By \cite{preprint:fol}, we have  $\classtop\simeq\Eqsheav{I_{\theory}}{M_{\theory}}$. Clearly, the forgetful functor $u^*:\Eqsheav{H_{\theory}}{M_{\theory}}\rightarrow \Eqsheav{I_{\theory}}{M_{\theory}}$ is faithful. The construction of Section \ref{Subsection: Equivariant sheaves on MT} can be seen to yield a conservative \theory-model $\cat{M}:\synt{C}{T}\rightarrow \Eqsheav{H_{\theory}}{M_{\theory}}$ also in the coherent case. Consequently, there is a geometric functor $m^*:\Eqsheav{I_{\theory}}{M_{\theory}}\rightarrow \Eqsheav{H_{\theory}}{M_{\theory}}$ classifying \cat{M}, and by inspection of the construction in \cite{preprint:fol}, $u^*\circ m^*\cong 1_{\Eqsheav{I_{\theory}}{M_{\theory}}}$ so that $u^*$ is essentially surjective on objects. Remains to show that $u^*$ is full. Let $\alg{R}=\pair{r:R\rightarrow M_{\theory},\rho}$, $\alg{S}=\pair{s:S\rightarrow M_{\theory},\sigma}$ be objects of \Eqsheav{H_{\theory}}{M_{\theory}}, and $F:u^*(\alg{R})\rightarrow u^*(\alg{S})$ a morphism of \Eqsheav{I_{\theory}}{M_{\theory}}. Let $\alg{M}\in M_{\theory}$, $x\in r^{-1}(\alg{M})$, and a homomorphism $\alg{h}:\alg{M}\rightarrow \alg{N}$ be given. Since every homomorphism in $H_{\theory}$ can be factored as an isomorphism followed by an inclusion, it can be assumed without loss that \alg{h} is an inclusion, $\alg{h}:\alg{M}\subseteq \alg{N}$. Then any neighborhood of \alg{M} must also contain \alg{N}. Now, the proof of Lemma \ref{Lemma: Well behaved section} goes through also for the coherent case, so we can choose a section $t:\bopen{\syntob{\alg{x}}{\phi},\alg{a}}\rightarrow R$  including $x$ such that for any $\alg{f}\in H_{\theory}$ such that $\alg{f}(\alg{a})=\alg{a}$ and with $d(\alg{f})\in \bopen{\syntob{\alg{x}}{\phi},\alg{a}}$ we have $\rho(\alg{f},t(d(\alg{f})))=t(c(\alg{f}))$. In particular, the section contains  $\rho(\alg{h},x)$. Applying $F$ to this section, we obtain a section of $S\rightarrow M_{\theory}$ containing $F(x)$ and $F(\rho(\alg{h},x))$. Using Lemma \ref{Lemma: Well behaved section} again, we find a section containing both $F(x)$ and $\sigma(\alg{h},F(x))$. Since $\alg{M}\subseteq \alg{N}$, the intersection of those two sections is non-empty over \alg{N}, so $F(\rho(\alg{h},x))=\sigma(\alg{h},F(x))$. Thus $F$ is also a morphism in \Eqsheav{H_{\theory}}{M_{\theory}}.
\end{proof}
\end{theorem}

\section*{Acknowledgements}
Thanks to Steve Awodey, Michael Makkai, Ji\v{r}\'{i} Rosick\'{y}, and Thomas Streicher for interesting and helpful discussions. This research was in part supported by the Eduard \v{C}ech Center through grant no.\ LC505.


\bibliographystyle{plain}
\bibliography{bibliografi}

\end{document}